\documentclass[12pt,reqno, oneside]{amsart}

\usepackage{latexsym}
\usepackage[centertags]{amsmath}
\usepackage{amsfonts}
\usepackage{amssymb}
\usepackage{amsthm}
\usepackage{newlfont}
\usepackage[square]{natbib}

\def\mytopsep{3mm}

\newtheoremstyle{myplain}{\mytopsep}{\mytopsep}{\itshape}{0pt}{\bfseries}{.}{3mm}{}
\newtheoremstyle{mydefinition}{\mytopsep}{\mytopsep}{\normalfont}{0pt}{\bfseries}{.}{3mm}{}

\newtheoremstyle{myremark}{\mytopsep}{\mytopsep}{\normalfont}{0pt}{\bfseries}{.}{3mm}{}

\theoremstyle{myplain}

\newtheorem{thm}{Theorem}[section]

\newtheorem{lem}[thm]{Lemma}

\theoremstyle{mydefinition}

\theoremstyle{myremark}

\allowdisplaybreaks[1]

\makeatletter\@addtoreset{equation}{section}\makeatother

\def\mb{\mathbf}

\def\NN{\mathbb{N}}

\def\CC{\mathbb{C}}

\def\ct{\mathop{\mathrm{CT}}}
\def\CT{\mathop{\mathrm{CT}}}

\renewcommand{\ll}{\langle\!\langle}
\renewcommand{\gg}{\rangle\!\rangle}
\newcommand{\lpp}{((}
\newcommand{\rpp}{))}
\newcommand{\lrq}[3]{\left(\frac{#1}{#2}#3\right)}

\newcommand{\qbinom}[2]{\genfrac{[}{]}{0pt}{}{#1}{#2}}
\newcommand{\qfac}[1]{(q)_{#1}}
\newcommand{\br}{\mathbf r}
\newcommand{\bk}{\mathbf k}
\newcommand{\Erk}{E_{\br,\bk}}
\newcommand{\Qrk}{\mathcal{Q}(b \Mid \br;\bk)}
\newcommand{\Mid}{\:|\:}  

\begin{document}

\title[A Short Proof of the Zeilberger-Bressoud $q$-Dyson
Theorem]{A Short Proof of the Zeilberger-Bressoud \\$q$-Dyson
Theorem}

\author{Ira M. Gessel$^*$}
\address{Department of Mathematics \\
Brandeis University\\
Waltham, MA 02454-9110} \email{gessel@brandeis.edu}

\author{Guoce Xin}
\address{Department
of Mathematics\\
Brandeis University\\
Waltham MA 02454-9110} \email{guoce.xin@gmail.com}

\date{February 22, 2005}
\thanks{$^*$Partially supported by NSF Grant DMS-0200596}
\begin{abstract}
We give a formal Laurent series  proof of Andrews's $q$-Dyson Conjecture, first proved by Zeilberger and Bressoud.
\end{abstract} \maketitle

\section{Introduction}

{}Freeman Dyson, motivated by a problem in particle physics,
conjectured the following identity in 1962:
\begin{thm}[Dyson's Conjecture]\label{t-dyson}
{}For nonnegative integers $a_0,a_1,\ldots ,a_n$,
\begin{equation}
\CT_{\mb{x}} \prod_{0\le i\ne j \le n}
\left(1-\frac{x_i}{x_j}\right)^{\!\!a_j} =
 \frac{(a_0+a_1+\cdots a_n)!}{a_0!\, a_1!\, \cdots a_n!},\label{e-dyson}
\end{equation}
where $\ct_{\mb{x}}$ denotes the constant term in $x_0,\dots,
x_n$.
\end{thm}

Dyson's conjecture was quickly proved by Wilson \cite{wilson} and
independently by Gunson~\cite{gunson}. An elegant recursive proof
was published by Good \cite{good1} in 1970.

A $q$-analog of Theorem \ref{t-dyson} was conjectured by George
Andrews \cite{andrews-qdyson} in 1975:
\begin{thm}[Zeilberger-Bressoud]\label{t-qdyson}
{}For any nonnegative integers $a_0,a_1,\dots ,a_n$,
\begin{align}\label{e-qdysonleft}
\CT_{\mb{x}} \prod_{0\le i<j\le n} \left({x_i\over
x_j}\right)_{\!\!a_i}\! \left({x_j\over x_i}q\right)_{\!\!a_j}
=\frac{(q)_{a_0+ \cdots+a_n}}{(q)_{a_0}\cdots
(q)_{a_n}},\end{align} where $(z)_m=(1-z)(1-qz) \cdots
(1-q^{m-1}z).$
\end{thm}

Andrews's $q$-Dyson conjecture attracted much interest, but was
not proved until 1985, in a combinatorial {\it tour de force\/} by
Zeilberger and Bressoud \cite{zeil-bressoud}. In related work,
Stanley \cite{stanley-qdyson1,stanley-qdyson2} reformulated the
conjecture in terms of symmetric functions and proved a limiting
form of the conjecture, and  Kadell \cite{kadell-qdyson} proved
the $4$-variable case using an approach similar to Good's.
Bressoud and Goulden \cite{breg} extended the Zeilberger-Bressoud
approach to some generalizations of  \eqref{e-qdysonleft}, and
Stembridge \cite{stembridge-qdyson} gave an elegant recursive
proof of the equal parameter case. Cherednik \cite{cherednik}
proved the Macdonald constant term conjecture for root systems
\cite{macdonald}, which generalizes the equal parameter case.

Zeilberger and Bressoud's combinatorial proof is, so far, the only proof of
Andrews's $q$-Dyson conjecture. We give here a very different and shorter
proof, using properties of formal Laurent series.

The idea behind the proof is the well-known fact that to prove the equality of
two polynomials of degree at most $d$, it is sufficient to prove that they are
equal at $d+1$ points. As is often the case, points at which the polynomials
vanish are most easily dealt with.

It is not difficult to show that for fixed nonnegative integers
$a_1,\dots, a_n$, both sides of \eqref{t-qdyson} are polynomials
in $q^{a_0}$ of degree at most $a_1+\cdots+a_n$ and that the
polynomial corresponding to the right side of \eqref{t-qdyson}
vanishes for $a_0= -1, -2,\dots, -(a_1+\cdots+a_n)$. The main part
of our proof is showing that the polynomial corresponding to the
left-hand side of \eqref{t-qdyson} also vanishes at these points;
we do this by expanding the left-hand side in partial fractions in
such a way that we can show that each summand has zero constant
term. The proof is completed by observing that the case $a_0=0$ of
\eqref{t-qdyson} is equivalent to the $n$-variable case.

\section{Basic Facts}
We use the following standard notation:
\begin{align}
\label{e-rfac}
(z)_n=\frac{(z)_\infty}{(zq^n)_\infty}= \prod_{m=0}^\infty
\frac{(1-zq^m)}{(1-zq^{m+n})}.
\end{align}

Note that if $p$ is a nonnegative integer, then
\begin{align*}
(z)_p &=(1-z)(1-zq)\cdots (1-zq^{p-1}) \\
(z)_{-p}&=\frac{1}{(1-zq^{-1})(1-zq^{-2})\cdots (1-zq^{-p})}.
\end{align*}

The $q$-binomial coefficients are defined for all integers $n$ and
nonnegative integers $m$ by
\begin{align}
\qbinom{n}{m} =\frac{(q^{n-m+1})_m}{(q)_m} =
\frac{(1-q^n)(1-q^{n-1})\cdots (1-q^{n-m+1})}{(1-q)(1-q^2)\cdots
(1-q^m)}.
\end{align}
{}For $0\le m\le n$, we have
\begin{align} \qbinom{n}{m}
=\frac{\qfac{n}}{\qfac{m}\qfac{n-m}}.
\end{align}

The well-known $q$-binomial theorem \cite[Theorem
2.1]{andrew-qbinomial} is the identity
\begin{align}
\label{e-qbinomial}
\frac{(az)_\infty}{(z)_\infty} = \sum_{k=0}^\infty \frac{(a)_k}{(q)_k} z^k.
\end{align}
Setting $z=uq^n$ and $a=q^{-n}$ in \eqref{e-qbinomial}, and using the definition
\ref{e-rfac}, we obtain
\begin{align}\label{e-qbinomialn}
(u)_n=\sum_{k=0}^\infty q^{k(k-1)/2}\qbinom{n}{k} (-u)^k
\end{align}
for all integers $n$.

We will also need the easily-proved identity
\begin{align}
\label{e-prod1} \left(\frac{x_i}{x_j}
\right)_{\!\!l}\!\left(\frac{x_j}{x_i}q
\right)_{\!\!m}=q^{\binom{m+1}{2}}\left(-\frac{x_j}{x_i}\right)^{\!\!m}
\!\left(\frac{x_i}{x_j} q^{-m}\right)_{\!\!l+m}.
\end{align}

\section{The Proof}
Let us fix $\mb{a}=(a_1,\dots ,a_n)$, where $a_1,\dots, a_n$
are nonnegative integers.   Let $a=a_1+\cdots +a_n$ and let
\begin{align}
P_{\mb{a}}(q^b)=\frac{\qfac{a}}{\qfac{a_1} \cdots \qfac{a_n}}
\qbinom{b+a}{a}=\frac{(1-q^{b+a})(1-q^{b+a-1})\cdots
(1-q^{b+1})}{\qfac{a_1} \cdots \qfac{a_n}}.
\end{align}
Then the right-hand side of \eqref{e-qdysonleft} is equal to
$P_{\mb{a}}(q^{a_0})$. We observe that
\begin{enumerate}
\item[(i)] $P_{\mb{a}}(q^b)$ is a polynomial in $q^b$ of degree at most $a$.
\item[(ii)] $P_{\mb{a}}(q^{b})=0$ for $q^b=q^{-1},q^{-2},\dots ,q^{-a}$.
\end{enumerate}
Moreover, $P_{\mb{a}}(q^b)$ is uniquely determined by these two
properties up to a constant factor (which may depend on $q$ but
not on $b$).

Let $Q_{\mb{a}}(q^{b})$ be defined by
\begin{align}\label{e-QAqb}
Q_{\mb{a}}(q^b)=\ct_{\mb{x}} \prod_{j=1}^n
\lrq{x_0}{x_j}{}_{\!\!b}\! \lrq{x_j}{x_0}{q}_{\!{}\!a_j}
\prod_{1\le i<j\le n} \lrq{x_i}{x_j}{ }_{\!\!a_i}
\!\lrq{x_j}{x_i}{q}_{\!\!a_j}.
\end{align}
Then the left-hand side of \eqref{e-qdysonleft} equals
$Q_{\mb{a}}(q^{a_0})$.

In fact $Q_{\mb{a}}(q^b)$ is well defined for negative integers
$b$ if we treat the rational function in \eqref{e-QAqb} as a
Laurent series in $x_0$. The following two lemmas show that
$Q_{\mb{a}}(q^b)$ equals $P_{\mb{a}}(q^ b)$ up to a constant
multiple.

\begin{lem}\label{l-main0}
{}For fixed $\mb{a}\in \NN^n$, $Q_{\mb{a}}(q^b)$ is a polynomial
in $q^b$ of degree at most $a$.
\end{lem}

\begin{lem}[Main Lemma]\label{l-main}
{}For any $\mb{a}\in \NN^n$, $Q_{\mb{a}}(q^b)=0$ for
$b=-1, -2, \dots ,-a$.
\end{lem}

Lemma \ref{l-main} is the heart of our proof of Theorem
\ref{t-qdyson}. We will prove it in Sections \ref{s-3} and
\ref{s-4}. We give here the rest of the proof of Theorem
\ref{t-qdyson}.

\begin{proof}[Proof of Lemma \ref{l-main0}]
By \eqref{e-prod1} we have
\begin{align*}
\left(\frac{x_0}{x_j} \right)_{\!\!b}\!\left(\frac{x_j}{x_0}q
\right)_{\!\!a_j}=q^{\binom{a_j+1}{2}}\left(-\frac{x_j}{x_0}\right)^{\!\!a_j}
\!\left(\frac{x_0}{x_j} q^{-a_j}\right)_{\!\!b+a_j},
\end{align*}
for all integers $b$, where both sides are regarded as Laurent series in $x_0$.

Thus $Q_{\mb{a}}(q^b)$ can be rewritten as
\begin{align}
\label{e-product} \ct_{\mb{x}}\: L(x_1,\dots ,x_n,\mb{a})
\prod_{j=1}^n q^{\binom{a_j+1}{2}}
\left(-\frac{x_j}{x_0}\right)^{\!\!a_j}\! \left(\frac{x_0}{x_j}
q^{-a_j}\right)_{\!\!b+a_j}\! ,
\end{align}
where $L(x_1,\dots,x_n,\mb{a})$ is a Laurent polynomial in $x_1,\dots,
x_n$ independent of $x_0$ and $b$.

Using the $q$-binomial theorem \eqref{e-qbinomialn}, we see that
for $1\le j\le n$,
$$ q^{\binom{a_j+1}{2}}
\left(-\frac{x_j}{x_0}\right)^{\!\!a_j} \left(\frac{x_0}{x_j}
  q^{-a_j}\right)_{\!\!b+a_j}
=\sum_{k_j\ge 0}C(k_j) \qbinom{b+a_j}{k_j}
x_0^{k_j-a_j}x_j^{a_j-k_j},$$ where
$C(k_j)=(-1)^{k_j}q^{\binom{a_j+1}2 + \binom {k_j}2 -k_ja_j}$.

Expanding the product in \eqref{e-product} and  taking the
constant term in $x_0$, we get
\begin{align} \label{e-midle} Q_{\mb{a}}(q^b)=\sum_{\mathbf{k}}
\qbinom{b+a_1}{k_1}\qbinom{b+a_2}{k_2}\cdots \qbinom{b+a_n}{k_n}
\ct_{x_1,\dots,x_n} L'(x_1,\dots ,x_n,\mb{a},\mathbf{k}),
\end{align}
where the sum ranges over all sequences $\mathbf{k}=(k_1,\dots,
k_n)$ of nonnegative integers such that $k_1+k_2+\cdots
+k_n=a_1+a_2+\cdots +a_n$, and $L'(x_1, \dots,
x_n,\mathbf{a},\mathbf{k})$ is a Laurent polynomial in $x_1,
\dots, x_n$ independent of $b$.
Since $\qbinom{b+a_i}{k_i}$ is a
polynomial in $q^b$ of degree $k_i$, each
summand in \eqref{e-midle} is a polynomial in $q^b$ of degree
$k_1+k_2+\cdots +k_n=a_1+a_2+\cdots +a_n$, and so is the sum.
\end{proof}

\begin{proof}[Proof of Theorem \ref{t-qdyson}]
We proceed by induction on $n$. Theorem \ref{t-qdyson} is trivial
for $n=0$ and reduces to the $q$-binomial theorem for $n=1$.
Suppose the theorem is true for $n$ variables. We may call these
variables $x_1,\dots, x_n$ rather than $x_0,\dots x_{n-1}$, so our induction
hypothesis implies that \eqref{e-qdysonleft} holds when $a_0=0$.

We will show that $P_{\mb{a}}(q^{a_0})=Q_{\mb{a}}(q^{a_0})$ for all
nonnegative integers $a_0$. We know that (i)
$P_{\mb{a}}(q^0)=Q_{\mb{a}}(q^0)$ by the induction hypothesis;
(ii) by Lemma \ref{l-main0}, both $P_{\mb{a}}(q^b)$ and
$Q_{\mb{a}}(q^b)$ define polynomials in $q^b$ of degree no greater
than $a$; (iii) by Lemma \ref{l-main},
$P_{\mb{a}}(q^b)=Q_{\mb{a}}(q^b)=0$ for $q^b=q^{-1},q^{-2},\dots,
q^{-a}$. So $P_{\mb{a}}(q^b)$ and $Q_{\mb{a}}(q^b)$ are equal as
polynomials in $q^b$.
\end{proof}

\section{Constant Term Evaluations}
\label{s-3}

We will evaluate the constant term
$Q_{\mb{a}}(q^b)$ defined by \eqref{e-QAqb}, where $b$ is a
negative integer, by partial fraction expansion. Although we are
taking the constant term of a Laurent series in $x_0$ with
coefficients that are Laurent polynomials in $x_1,\dots, x_n$,
when we expand by partial fractions, we get terms that are not of
this form, and in order to evaluate their constant terms we need
to work in a larger ring: the field of iterated Laurent series
$K\ll x_n,x_{n-1},\dots ,x_0 \gg=K\lpp x_n\rpp\lpp
x_{n-1}\rpp\dots\lpp x_0\rpp$, where $K=\CC(q)$, in which all
series are  regarded first as Laurent series in $x_0$, then as
Laurent series in $x_1$, and so on. For a more detailed account of
the properties of this field, with other applications, see
\cite{xinresidue} and \cite{xiniterate}.

Every element of  $K\ll x_n,x_{n-1},\dots ,x_0 \gg$ has a unique
Laurent series expansion. The series expansions of $1/(1-q^k
x_i/x_j)$ will be especially important. If $i<j$ then
$$\frac{1}{1-q^k x_i/x_j}=\sum_{l= 0}^\infty q^{kl} x_i^l x_j^{-l}.$$
However, if $i>j$ then this expansion is not valid and instead we
have the expansion
$$ \frac{1}{1-q^k x_i/x_j}=\frac1{-q^k x_i/x_j(1-q^{-k}x_j/x_i)}
    =\sum_{l=0}^\infty -q^{-k(l+1)} x_i^{-l-1}x_j^{l+1}.$$

{}For $F(\mathbf{x}) \in K\ll x_n,x_{n-1},\dots ,x_0 \gg$, the
\emph{constant term} of $F(\mathbf{x})$ in $x_i$, denoted by
$\ct_{x_i} F(\mathbf{x})$, is defined to be the sum of those terms
in the series expansion of $F(\mathbf{x})$ that are free of $x_i$.
This definition clearly extends the constant term operators used
earlier. It follows that
\begin{equation}
\label{e-ct} \ct_{x_i} \frac{1}{1-q^k x_i/x_j} =
\begin{cases}
    1, & \text{ if }i<j, \\
    0, & \text{ if }i>j. \\
\end{cases}
\end{equation}
We shall call the monomial $M=q^k x_i/x_j$ \emph{small} if $i<j$
and \emph{large} if $i>j$.  Thus the constant term in $x_i$ of
$1/(1-M)$ is 1 if $M$ is small and $0$ if $M$ is large.

An important property of the constant term operators defined in
this way is their commutativity:
$$\ct_{x_i} \ct _{x_j} F(\mathbf{x}) = \ct_{x_j} \ct _{x_i} F(\mathbf{x}).$$
Commutativity implies that the constant term in a set of variables
is well-defined, and this property will be used in our proof of
the Main Lemma. (Note that, by contrast, the constant term
operators in \cite{zeil} do not commute.)

The \emph{degree} of a rational function of $x$ is the degree in
$x$ of the numerator minus the degree  in $x$ of the denominator.
{}For example, if $i\ne j$ then  the degree  of
$1-x_j/x_i=(x_i-x_j)/x_i$ is $0$ in $x_i$ and $1$ in $x_j$. A
rational function is called \emph{proper} in $x$ if its degree in
$x$ is negative. The following lemma gives a formula for the
constant term in $x_k$ of certain elements of $K\ll
x_n,x_{n-1},\dots ,x_0 \gg$ which are proper rational functions of
$x_k$.

\begin{lem}\label{l-vvv}
Let
$$R=\frac{p(x_k)}{x_k ^d \prod_{i=1}^m (1-x_k/\alpha_i)},$$
be a proper rational function of $x_k$ where $p(x_k)$ is a
polynomial in $x_k$, and the $\alpha_i$ are distinct monomials,
each of the form $x_t q^s$. Then
\begin{align}
\ct_{x_k} R=\sum_j  \bigl(R\, (1-x_k/\alpha_j)\bigr)\Bigr|_{x_k
=\alpha_j},
\end{align} where the sum ranges
over all $j$ such that $x_k/\alpha_j$ is small.
\end{lem}
\begin{proof}
The field $K\ll x_n,\dots ,x_0\gg $ contains the polynomial ring
$K[x_0,\dots ,x_n]$ as a subring and hence contains the field
$K(x_0,\dots ,x_n)$ of rational functions as a subfield. Thus
any identity in $K(x_0,\dots ,x_n)$ is also an identity in $K\ll x_n,\dots
,x_0\gg$.

The partial fraction decomposition
of $R$ with respect to $x_k$ is
$$\frac{p(x_k)}{x_k^d \prod_{i=1}^m (1-
x_k/\alpha_i)}=\frac{p_0(x_k)}{x_k^d}+\sum_{j=1}^m \frac{1}{1-
x_k/\alpha_j}  \left. \left(\frac{p(x_k)}{x_k^d \prod_{i=1,i\ne
j}^m (1-x_k/\alpha_i)}\right)\right|_{x_k=\alpha_j},$$ where
$p_0(x_k)$ is a polynomial in $x_k$ of degree less than $d$. The
term $p_0(x_k)/x_k^d$ contributes nothing to the constant term in
$x_k$, and $1/(1-x_k/\alpha_j)$ contributes to the constant term
in $x_k$ only if $x_k/\alpha_j$ is small. The result of the lemma
then follows easily.
\end{proof}

The following lemma plays an important role in our argument.
\begin{lem}\label{l-tournament}
Let $A_1,\dots, A_s$ be nonnegative integers. Then for any
positive integers $k_1,\dots, k_s$ with $1\le k_i\le A_1+\cdots
+A_s$ for all $i$, either $1\le k_i\le A_i$ for some $i$ or
$-A_j\le k_i-k_j\le A_i-1$ for some $i<j$.
\end{lem}
\begin{proof}
We prove by contradiction that there is no $k_1,\dots ,k_s$ such
that for all $i$, $A_i < k_i\le A_1+\cdots +A_s$, and for all
$i<j,$ either $k_i-k_j\ge A_i$ or $k_i-k_j\le -A_j-1$. Suppose
$k_1,\dots ,k_s$ do satisfy these conditions. We construct a
tournament on $1,2,\dots ,s$ with numbers on the arcs as follows:
{}For $i<j$, if $k_i-k_j\ge A_i$ then we draw an arc
$i\mathop{\longleftarrow}\limits^{A_i} j$ from $j$ to $i$ and if
$k_i-k_j\le -1-A_j$ then we draw an arc
$i\mathop{\longrightarrow}\limits^{A_j+1} j$ from $i$ to $j$.

We call an arc from $u$ to $v$ an \emph{ascending arc} if $u<v$
and a \emph{descending arc} if $u>v$.  We note two facts: (i) the
number on an arc from $u$ to $v$ is less than or equal to
$k_v-k_u$, and (ii) the number on an ascending arc is always
positive.

A  consequence of (i) is that for any directed path from $e$ to
$f$, the sum along the arcs is less than or equal to $k_f-k_e$. It
follows that the sum along a cycle is nonpositive. But any cycle
must have at least one ascending arc, and by (ii) the number on
this arc is positive, and so the sum along the cycle is positive.
Thus there can be no cycles.

Therefore the tournament we have constructed is transitive, and
hence defines a total ordering $\rightarrow$ on $1,2,\dots ,s$.
Assume the total ordering is given by $i_1\rightarrow
i_2\rightarrow \cdots \rightarrow i_{s-1}\rightarrow i_s$. Then
$k_{i_s}-k_{i_1}\ge A_{i_2}+A_{i_3}+\cdots +A_{i_s}$. This implies
that
\begin{align*}
k_{i_s}&\ge k_{i_1}+A_{i_2}+A_{i_3}+\cdots+A_{i_s}\\
    & > A_{i_1}+A_{i_2}+A_{i_3}+\cdots +A_{i_s}\\
    &=A_1+A_2+\cdots +A_s,
\end{align*}
 a contradiction.
\end{proof}

\section{Proof of the Main Lemma}
\label{s-4}

Let
\begin{align*}
\mathcal{Q}(b)=\prod_{j=1}^n \frac{(x_jq/x_0)_{a_j}} {(1-x_0/x_jq)
(1-x_0/x_jq^2) \cdots (1-x_0/x_jq^b)} \prod_{1\le i <j \le n}
\lrq{x_i}{x_j}{}_{\!\!a_i}\!\lrq{x_j}{x_i}{q}_{\!\!a_j},
\end{align*}
so that in the notation of the previous section, $\ct_{\mb{x}}\:
\mathcal{Q}(b)=Q_{\mb{a}}(q^{-b})$. Our goal is to show that
$\ct_{\mb{x}}\mathcal{Q}(b)=0$ for $b=1,2,\dots, a$.

Since the degree in $x_0$ of $1-x_jq^i/x_0=(x_0-x_j q^i)/x_0$ is
$0$, $\mathcal{Q}(b)$ is proper  in $x_0$, with degree $-nb$.

Applying Lemma \ref{l-vvv}, we have
\begin{equation}
\label{e-P1} \ct_{x_0}\mathcal{Q}(b)=\sum_{0<r_1\le n\atop 1\le
k_1\le b} \mathcal{Q}(b\Mid r_1;k_1),
\end{equation}
where
\begin{equation*}\label{e-QQ1}
\mathcal{Q}(b \Mid r_1;k_1)= \left. \mathcal{Q}(b)
\left(1-\displaystyle\frac{x_0}{x_{r_1}
q^{k_1}}\right)\right|_{x_0=x_{r_1}q^{k_1}}.
\end{equation*}

For each term in \eqref{e-P1} we will extract the constant term in
$x_{r_1}$, and then perform further constant term extractions,
eliminating one variable at each step. In order to keep track of
the terms we obtain, we introduce some notation.

For any rational function $F$ of $x_0, x_1,\dots, x_n$, and for
sequences of integers $\bk=(k_1,k_2,\dots, k_s)$ and $\br=(r_1,
r_2,\dots, r_s)$ let $E_{\br,\bk}\,F$ be the result of replacing
$x_{r_i}$ in $F$ with $x_{r_s}q^{k_s -k_i}$ for $i=0,1,\dots,
s-1$, where we set $r_0=k_0=0$. Then for $0<r_1<r_2<\cdots <r_s\le
n$ and $0<k_{i}\le b$,  we define
\begin{equation}
\label{e-QQ2} \Qrk=\mathcal{Q}(b \Mid r_1,\dots, r_s; k_1, \dots,
k_s)=\Erk\, \left[\mathcal{Q}(b)
 \prod_{i=1}^s \left(1-\displaystyle\frac{x_0}{x_{r_i}
q^{k_i}}\right)\right].
\end{equation}
Note that the product on the right-hand side of \eqref{e-QQ2}
cancels all the factors in the denominator of $\mathcal{Q}$ that
would be taken to zero by $\Erk$.

\begin{lem}
The rational functions $\Qrk$ have the following two properties:
\label{l-Qrk}
\begin{enumerate}
\item[(i)] If $1\le k_{i}\le a_{r_1}+\cdots +a_{r_s}$ for
     all $i$ with $1\le i\le s$, then
    $\Qrk=0$.
\item[(ii)] If $k_{i}>a_{r_1}+\cdots +a_{r_s}$ for some $i$ with $1\le i\le s$
     and     $n>s$, then
\begin{equation}
\label{e-Qii} \ct_{x_s} \Qrk =\sum_{r_s < r_{s+1}\le n\atop 1\le
k_{{s+1}}\le b}
    \mathcal{Q}(b\Mid r_1,\dots, r_s, r_{s+1};k_1,\dots, k_s,k_{s+1}).
\end{equation}
    \end{enumerate}
\end{lem}

\begin{proof}[Proof of property (i)]
By Lemma \ref{l-tournament}, either  $1\le k_i\le a_{r_i}$ for
some $i$ with $1\le i \le s$, or $-a_{r_j}\le k_i-k_j\le
a_{r_i}-1$ for some $i<j$. If $1\le k_i\le a_{r_i}$ then $\Qrk$
has the factor
$$\Erk
\left[\lrq{x_{r_i}}{x_{0}}{q}_{\!\!a_{r_i}} \right]
=\lrq{x_{r_s}q^{k_s-k_i}}{x_{r_s}q^{k_s}}{q}_{\!\!a_{r_i}}
=(q^{1-k_i})_{a_{r_i}}=0.$$

If $-a_{r_j}\le k_i-k_j\le a_{r_i}-1$ where $i<j$ then $\Qrk$ has
the factor
$$\Erk\, \left[\lrq{x_{r_i}}{x_{r_j}}{}_{\!\!a_{r_i}}\!\!\lrq{x_{r_j}}{x_{r_i}}{q}_{\!\!a_{r_j}}\right], $$
which by \eqref{e-prod1} is equal to
$$q^{\binom{a_{r_j}+1}{2}} \left(-\frac{x_{r_j}}{x_{r_i}}\right)^{a_{r_j}}\!
\lrq{x_{r_i}}{x_{r_j}}{q^{-a_{r_j}}}_{\!\!a_{r_i}+a_{r_j}}
\!\!\!\!=q^{\binom{a_{r_j}+1}{2}}
(-q^{k_i-k_j})^{a_{r_j}}(q^{k_j-k_i-a_{r_j}})_{a_{r_i}+a_{r_j}}=0.
$$

\smallskip
\noindent\emph{Proof of property (ii).} Note that since $b\ge k_i$
for all $i$, the hypothesis implies that   $b>a_{r_1}+\cdots
+a_{r_s}$.

We first show that $\mathcal{Q}(b\Mid \mb{r};\mb{k})$
is proper in $x_{r_s}$. To do this we write
$\Qrk$ as $N/D$, in which $N$ (the ``numerator") is
$$\Erk\, \Biggl[\prod_{j=1}^n\lrq{x_j}{x_0}q_{\!\!a_j}
\cdot \prod_{\substack{1\le i, j\le n\\ j\neq i}}
\left(\frac{x_i}{x_j}\,q^{\chi(i>j)}\right)_{\!\!a_i}\Biggr],$$
where $\chi(S)$ is $1$ if the statement $S$ is true, and $0$
otherwise, and $D$ (the ``denominator") is
$$\Erk\, \left[\prod_{j=1}^n \lrq{x_0}{x_jq^b}{}_b\biggm/\prod_{i=1}^s\left(1-\frac
{x_0}{x_{r_i}q^{k_i}}\right)\right].$$ Now let
$R=\{r_0,r_1,\dots,r_s\}$. Then the degree in $x_s$ of
$$\Erk\, \left[\left(1-\frac{x_i}{x_j}q^l\right)\right]$$
is 1 if $i\in R$ and $j\not\in R$, and is 0 otherwise, as is
easily seen by checking the four cases. Thus the  part of $N$
contributing to the degree in $x_{r_s}$ is
$$\Erk\,\left[\prod_{i=1}^s \prod_{j\ne r_0,\dots ,r_s}
\left(\frac{x_{r_i}}{x_j}q^{\chi(r_i>j)}\right)_{\!\!a_{r_i}}\right],$$
 which has degree $(n-s)(a_{r_1}+\cdots +a_{r_s})$, and the
part of $D$ contributing to the degree in $x_{r_s}$ is
$$\Erk\,\biggl[\,\prod_{j\ne r_0,\dots,r_s}\lrq{x_0}
{x_jq^b}{}_{\!\!b}\biggr],$$ which has degree $(n-s)b$.

Thus the total degree of $\Qrk$ in $x_{r_s}$ is
$(n-s)(a_{r_1}+\cdots +a_{r_s} - b)<0$, so $\Qrk$ is proper in
$x_{r_s}$.

Next we apply Lemma \ref{l-vvv}. For any rational function $F$ of
$x_{r_s}$ and integers $j$ and $k$,  let $T_{j,k} F$ be the result
of replacing $x_{r_s}$ with $x_{j}q^{k-k_s}$ in $F$. Since
$x_{r_s}q^{k_s}/(x_jq^k)$ is small when $j>r_s$ and is large when
$j<r_s$, Lemma \ref{l-vvv} gives
\begin{equation}
\label{e-TQ} \ct_{x_s} \Qrk =\sum_{r_s < r_{s+1}\le n\atop 1\le
   k_{{s+1}}\le b} T_{r_{s+1},k_{s+1}} \left[\Qrk
   \left(1-\frac{x_{r_s}q^{k_s}}{x_{r_{s+1}}q^{k_{s+1}}}\right)\right].
\end{equation}
We must show that the right-hand side of $\eqref{e-TQ}$ is equal
to the right-hand side of \eqref{e-Qii}. Let us set
$\br'=(r_1,\dots, r_s, r_{s+1})$ and $\bk'=(k_1,\dots, k_s,
k_{s+1})$. Then the equality follows easily from the identity
\begin{equation}
\label{e-TE} T_{r_{s+1},k_{s+1}}\circ \Erk= E_{\br',\bk'}.
\end{equation}
To see that \eqref{e-TE} holds, we have
$$
(T_{r_{s+1},k_{s+1}}\circ \Erk)\, x_{r_i}
  =T_{r_{s+1},k_{s+1}}\, \left[ x_{r_s}q^{k_s-k_i}\right]
  = x_{r_{s+1}}q^{k_{s+1}-k_i}= E_{\br',\bk'}\,  x_{r_i},
$$
and if $j\notin\{r_0,\dots, r_s\}$ then $(T_{r_{s+1},k_{s+1}}\circ
\Erk)\, x_{j}=x_j=  E_{\br',\bk'}\, x_{j}$.

\end{proof}

\vspace{3mm}
\begin{proof}[Proof of the Main Lemma]
We prove by induction on $n-s$ that
$$\ct_{\mb{x}} \Qrk = 0;$$
the Main Lemma is the case $s=0$. (Note that taking the constant
term with respect to a variable that does not appear has no
effect.)

We may assume that $s\le n$ and $0<r_1<\cdots<r_s\le n$, since
otherwise $\Qrk$ is not defined. If $s=n$ then $r_i$ must equal
$i$ for $i=1,\dots ,n$ and thus $\Qrk=\mathcal{Q}(b\Mid 1,2,\dots,
n; k_{1}, k_{2},\dots, k_{n})$, which is 0 by part (i) of
Lemma \ref{l-Qrk}, since $k_i\le b\le a_1+\cdots + a_n$ for each
$i$.

Now suppose that $0\le s<n$. If part (i) of Lemma \ref{l-Qrk}
applies, then $\Qrk=0$. Otherwise, part (ii) of Lemma \ref{l-Qrk}
applies and \eqref{e-Qii} holds. Applying $\ct_{\mb{x}}$ to both
sides of \eqref{e-Qii} gives
$$
\ct_{\mb{x}}\Qrk
  =\sum_{r_s < r_{s+1}\le n\atop 1\le k_{{s+1}}\le b}
   \ct_{\mb{x}}\mathcal{Q}
   (b\Mid r_1,\dots, r_s, r_{s+1};k_1,\dots,k_s,k_{s+1}),
$$
and by induction, every term on the right is 0.

\end{proof}

\bibliographystyle{amsplain}

\providecommand{\bysame}{\leavevmode\hbox
to3em{\hrulefill}\thinspace}
\providecommand{\MR}{\relax\ifhmode\unskip\space\fi MR }
\providecommand{\MRhref}[2]{%
  \href{http://www.ams.org/mathscinet-getitem?mr=#1}{#2}
} \providecommand{\href}[2]{#2}

\end{document}